\documentclass{amsart}
\usepackage{amsmath}
\usepackage{amsfonts}
\usepackage{amssymb}
\usepackage{graphicx}
\usepackage{multirow}
\usepackage{stix}
\usepackage{courier}

\numberwithin{equation}{section}
\begin{document}

\title[Isogonic and isodynamic points of a simplex in a real affine space ]{Isogonic and isodynamic points of a simplex \protect\\ in a real affine space}
\author{Manfred Evers}
\curraddr[Manfred Evers]{Bendenkamp 21, 40880 Ratingen, Germany}
\email[Manfred Evers]{manfred\_evers@yahoo.com}
\date{\today}

\begin{abstract}A non-equilateral triangle in a Euclidean plane has exactly two isogonic and two isodynamic points. There are a number of different but equivalent characterizations of these triangle centers. The aim of this paper is to work out characteristic properties of isogonic and isodynamic centers of simplices that can be transferred to higher dimensions. In addition, a geometric description of the Weiszfeld algorithm for calculating the Fermat point of a simplex is given.
\end{abstract}\vspace*{1 mm}

\maketitle

\noindent\textbf{1. Introduction and terminology.}\hspace*{\fill} \vspace*{1 mm}\\
\noindent\hspace*{5mm}Let $A_1,\dots,A_{n+1}$ be $\,n+1$ affinely independent points in a real Euclidean-affine space $\mathcal{A}$ of dimension $n>1$. We work with barycentric coordinates with respect to the tupel $(A_1,\dots,A_{n+1})$. Given a point $P\in \mathcal{A}$, we write $P = [p_1,\dots,p_{n+1}]$ \,(resp. $P = [p_1:\cdots:p_{n+1}]\,)$, if $p_1,\dots,p_{n+1}$ are the absolute \,(\,resp. homogeneous\,)\, coordinates of $P$ with respect to $(A_1,\dots,A_{n+1})$.
The set $\Sigma = \bigl\{P=[p_1,\dots,p_{n+1}]\; |\; p_i \ge 0\; \textrm{for}\; i=1,\dots,n{+}1\bigr\}$ is an $n$-simplex. The affine subspace $A_{i_1}\sqcup\dots\sqcup A_{i_{k+1}}$ spanned by $k{+}1$ points $A_{i_1},\dots,A_{i_{k+1}} \in \,\{A_1,\dots,A_{n+1}\}$ is called $k$-\textit{sideplane} of $\Sigma$, whilst the $k$-simplex $\Sigma \cap (A_{i_1}\sqcup\dots\sqcup A_{i_{k+1}}) =: \Sigma_{i_1,\dots,\,i_{k+1}}$ is called a $k$-face of $\Sigma$. Instead of $1$-sideplanes we usually speak of \textit{sidelines}, and the $1$-faces are also called \textit{edges}, the $(n{-}1)$-faces  \textit{facets} of $\Sigma$.\\
Let $d_{ij}$ be the distance between two vertices $A_i,A_j$ of $\Sigma$, then the squared distance between two points $P = [p_1,\dots,p_{n+1}]$ and $Q = [q_1,\dots,q_{n+1}]$ is given by \vspace*{2mm}\\
\centerline{$\displaystyle d^2(P,Q) = \;-\!\sum_{{1\,\le\,i\,<\,j\,\le\,n{+}1}} (d_{ij})^2(p_i-q_i)(p_j-q_j)\hspace*{10mm}(\star)\; .$}\\
\textit{Remarks}.\\
$\bullet\;$ If $\mathcal{A}$ is the real vector space $\mathbb{R}^{n}$ with Euclidean dot product $\cdot$ and Euclidean distance function $d(P,Q)=  \sqrt{(P-Q)\cdot(P-Q)}\,$, then the equation $(\star)$ applies, as it was shown by Coxeter [2].\\
$\bullet\;$ As we know from Galilean and from Lorentz-Minkowski geometry, squared distances between points can take values that are non-positive real numbers (cf. [8]). But in Euclidean-affine spaces squared distances have to be positive, and, moreover, volumes of simplices of dimension $\le n$ take only positive values. \\

\noindent\textbf{2. Properties of the isogonic and of the isodynamic points in the planar case n=2.}\hspace*{\fill}\vspace*{1 mm} \\
\noindent\hspace*{5mm} In this paper we want to find properties that characterize isogonic and isodynamic points of an $n$-simplex. At first, we examine the $2$-dimensional case. Instead of giving a definition for these centers, we list several well known properties which characterize these points in a unique way, cf. [3,10,13,17,19,23].\vspace*{-3.5 mm}\\

A point $P$ is an \textit{isodynamic point} of a triangle $A_1A_2A_3$ if any of the following (equivalent) statements is true:\\
(1)  The mirror image of triangle $A_1A_2A_3$ with respect to a circle with center $P$ is an equilateral triangle. \\
(2)  Let $C_i$ denote the intersection of the tripolar of the symmedian $K$ ($K$ defined as the isogonal conjugate of the centroid $G$) with the sideline opposite the vertex $A_i$ and let $\mathscr{C}_i$ denote the circle with center $C_i$ passing through the vertex $A_i$, $i\,=\,1,2,3$; then $P$ is a common point of the three circles $\mathscr{C}_1, \mathscr{C}_2, \mathscr{C}_3.$\\
(3)  Let $\mathscr{C}$ denote the circumcircle of triangle $A_1A_2A_3$ and $\textrm{T}_{A_i}\!\mathscr{C}$ denote the tangent of $\mathscr{C}$ at the point $A_i$, let $C_{\;i}^\star$ denote the intersection of $\textrm{T}_{A_i}\! \mathscr{C}$ with the sideline opposite $A_i$ and $\mathscr{C}_{\;i}^\star$ denote the circle with center $C_{\;i}^\star$ passing through the vertex $A_i$, $i\,=\,1,2,3$; then $P$ is a common point of the three circles $\mathscr{C}_{\;1}^\star, \mathscr{C}_{\;2}^\star, \mathscr{C}_{\;3}^\star.$\\
(4)  The pedal triangle of $P$ is equilateral.\\
(5) $P$ is a point that satisfies the equation $d(P,A_1)\,d_{23} = d(P,A_2)\,d_{13} = d(P,A_3)\,d_{12}\,$.\vspace*{-2mm}\\

Let us assume that the triangle $A_1A_2A_3$ is not equilateral. Then there exist precisely two isodynamic points $J_1, J_2$, and the following statements hold:\vspace*{-2mm}\\

\noindent(6) $J_1, J_2$ both lie on the Brocard axis, a line through the symmedian $K$ and the circumcenter $O$ of triangle $A_1A_2A_3$.\\
(7) One of the two isodynamic points is a point inside the triangle, and the other is the mirror image of this point with respect to the circumcircle.\\
(8) The  points $O,J_1,K,J_2$ form a harmonic range.\vspace*{-2mm}\\

The circles $\mathscr{C}_{\;1}^\star, \mathscr{C}_{\;2}^\star, \mathscr{C}_{\;3}^\star$ are called the Apollonian circles of the triangle $A_1A_2A_3$. A point $Q$ is a point on $\mathscr{C}_{\;i}^\star$ precisely when\vspace*{1mm}\\
\centerline{$\displaystyle \frac{d(Q,A_j)}{d(Q,A_k)} = \frac{d_{ij}}{d_{ik}},\; \{i,j,k\} = \{1,2,3\}\;\;\;\;\;\;\;\;(\star\star)\,.$}\\

\noindent Let $\textrm{T}_{J_k}\!\mathscr{C}_{\;i}^\star$ be the tangent of $\mathscr{C}_{\;i}^\star$ at the point $J_k$, $k = 1,2$ and $i=1,2,3$. Then:\vspace*{-2mm}\\

\noindent(9) The three singular conics $\textrm{T}_{J_k}\!\mathscr{C}_{\;1}^\star\,\cup\,\textrm{T}_{J_k}\!\mathscr{C}_{\;2}^\star, \textrm{T}_{J_k}\!\mathscr{C}_{\;2}^\star\,\cup\,\textrm{T}_{J_k}\!\mathscr{C}_{\;3}^\star, \textrm{T}_{J_k}\!\mathscr{C}_{\;3}^\star\,\cup\,\textrm{T}_{J_k}\!\mathscr{C}_{\;3}^\star$ are congruent in pairs, k = 1,\,2. In other words, these are unions of two lines with the same angle of intersection.\\

A point $P$ is an \textit{isogonic point} (often also called \textit{isogonic center}) of triangle $A_1A_2A_3$ if any of the three following (equivalent) statements is true:\\
(10)  The three singular conics $(P\sqcup A_1)\cup (P\sqcup A_2), (P\sqcup A_2)\cup (P\sqcup A_3), (P\sqcup A_3)\cup (P\sqcup A_1)$ are congruent in pairs.\\
(11)  The inversion of triangle $A_1A_2A_3$ in a circle with center $P$ leads to an equilateral triangle.\\
(12)  The antipedal triangle of $P$ is equilateral.\\

There exist exactly two isogonic points, $F_1$ and $F_2$. These are the isogonal conjugates of the points $J_1$ and $J_2$, respectively.\\
If $F_1$ is inside the triangle $A_1A_2A_3$, then it is the Fermat-Torricelli point of the triangle, it minimizes the function $P\mapsto d(P,A_1)+d(P,A_2)+d(P,A_3)$. \\

\noindent\textit{Remarks}. Let us look at the situation in elliptic and in hyperbolic planes. In these planes there are two points - let us call them $J_1$ and $J_2$ - which satisfy (2). Both are points on the line $K\sqcup O$.  One more point on this line is the Lemoine point $\tilde{K}$, a point whose tripolar is orthogonal to $K\sqcup O$ and meets the sidelines of triangle $A_1A_2A_3$ in $C_{\;1}^\star, C_{\;2}^\star, C_{\;3}^\star$. In the affine case, the symmedian agrees with the Lemoine point, but in planes with nonzero Gaussian curvature these are two different triangle centers, see [9].\\
The two points $J_1$ and $J_2$ also satisfy (9), and the four points $O,J_1,\tilde{K},J_2$ form a harmonic range. 
But statements (1), (3), (4), (7) do not apply, in general. The equation given in (5) is responsible for the name of the two centers; it has to be replaced in the hyperbolic case by the equation\\
\centerline{$\sinh(\frac{1}{2}d(P,A_1))\sinh(\frac{1}{2}d_{23}) = \sinh(\frac{1}{2}d(P,A_2))\sinh(\frac{1}{2}d_{13}) = \sinh(\frac{1}{2}d(P,A_3))\sinh(\frac{1}{2}d_{12}) $}\\
and in the elliptic case by\\ 
\centerline{$\sin(\frac{1}{2}d(P,A_1))=\sin(\frac{1}{2}d_{23}) = \sin(\frac{1}{2}d(P,A_2))\sin(\frac{1}{2}d_{13}) = \sin(\frac{1}{2}d(P,A_3))\sin(\frac{1}{2}d_{12}) $\,.}\\
The name \textit{Apollonian circles} for the circles $\mathscr{C}_{\;1}^\star, \mathscr{C}_{\;2}^\star, \mathscr{C}_{\;3}^\star$ is justified also in planes with nonzero Gaussian curvature, because their equations are very similar to equation $(\star\star)$, cf.\,[9].\vspace*{0.5mm}
 
In  elliptic and in  hyperbolic planes, there are two points, $F_1$ and $F_2$, say, for which statement (10) is true (cf.\,[9]), and if one of these points $F_1, F_2$ lies inside the triangle, this point is the Fermat-Torricelli point of the triangle (cf. [11]). There are strong indications (based on experiments with GeoGebra) that $F_1$ and $F_2$ are exactly the points that satisfy condition (11), see [9].  But (12) does not necessarily apply to them, and $F_1$ and $F_2$ are, in general, not isogonal conjugates of $J_1$ and $J_2$. \vspace*{-1mm}\\

What is the situation like in a Lorentz-Minkowski plane or in a Galilean plane? 
Apollonian circles exist in a Lorentz-Minkowski plane, however their common points are points on the line at infinity. In a Galilean plane, each Apollonian circle consists of two parallel lines and all these six lines meet at the absolute pole. (A nonsingular circle touches the line at infinity at the absolute pole, and this point is the center of the circle.) So there are no isodynamic points in Galilean and Lorentz-Minkowski planes, nor are there isogonic points. \\

\noindent\textbf{3. Generalized Apollonian spheres.}\hspace*{\fill}\vspace*{1 mm} \\
In this section we generalize results and ideas published by P. Yiu [24] to higher dimensions.\vspace*{-2mm}\\

We come back to the general case of an $n$-simplex $\Sigma$ with vertices $A_1,\dots A_{n+1}$. Let $P = [p_1{:}\cdots{:}p_{n+1}]$ be a point not on any of the $(n{-}1)$-sideplanes of $\Sigma$, thus $p_1p_2\cdots p_{n+1} \ne 0$. The $(n-1)$-plane $P\sqcup A_3\sqcup \dots \sqcup A_{n+1}$ meets the line $A_1 \sqcup A_2$ at the point $P_{12}:=
[p_1:p_2:0:\cdots:0]$, while the $\Sigma$-polar plane of $P$, this is the $(n{-}1)$-plane $P^\Sigma := \big{\{}[x_1{:}\cdots {:}x_{n+1}]\;|\;{\frac{x_1}{p_1}+\dots+\frac{x_{n+1}}{p_{n+1}}=0} \big{\}}$, meets the line $A_1 \sqcup A_2$ at the point $P_{12}^\star:= [{-p_1{:}p_2{:}}\\ {0:\cdots:0]}$. The midpoint of these two points, $Q_{12}:= \frac{1}{2}P_{12}+\frac{1}{2}P_{12}^\star = {[-p_1^2{:}p_2^2{:}0{:}\cdots{:}0]}$, is a point on the $\Sigma$-polar plane of the barycentric square $P^2$ of $P$, as can be easily checked. The $(n{-}1)$-sphere with diameter $[\!\textbf{[} P_{12},P_{12}^\star {\textbf{]}\!]} := \big{\{}Q={t P_{12}+(1{-}t)P_{12}^\star}\;| \;0\le t\le 1\big{\}}$ and center $Q_{12}$ is denoted by $\mathcal{S}_{12}$.\vspace*{1mm}\\
This sphere $\mathcal{S}_{12}$ meets the circumsphere of $\Sigma$ orthogonally.\\
\textit{Proof}. The points $A_1, P_{12}, A_2, P_{12}^\star$ form a harmonic range. Therefore the points $A_1, A_2$ are inversive with respect to $\mathcal{S}_{12}$, and every sphere through $A_1$ and $A_2$, especially the circumsphere of $\Sigma$, is orthogonal to $\mathcal{S}_{12}.\;\,$  $\Box$\vspace*{1mm}

The points $P_{ij},P_{ij}^\star,Q_{ij}$ and $(n{-}1)$-spheres $\mathcal{S}_{ij}$, $1\le i\le j\le n{+}1$  are defined likewise. Following Yiu, we will call these spheres \textit{generalized Apollonian spheres} of $\Sigma$ with respect to $P$. The sphere $\mathcal{S}_{ij}$ is the locus of points  $R$ satisfying $d(A_i, R) : d(A_j;R) = 1/|p_i| : 1/|p_j|.$ From this follows that $R$ is a common point of two spheres $\mathcal{S}_{ij}, \mathcal{S}_{jk}$, $1\le i<j<k\le n+1$, precisely when $d(A_i, R) : d(A_j; R) : d(A_k; R) = 1/|p_i| : 1/|p_j| : 1/|p_k|$. But then $R$ must be a point on $\mathcal{S}_{ik}$, as well. The radical $(n{-}1)$-plane of the three spheres is perpendicular to the line through their centers and contains the circumcenter $O$. Therefore, the intersection of all these radical $(n{-}1)$-planes is the line through $O$ perpendicular to the $\Sigma$-polar plane of $P^2$. As a consequence, all $(n{-}1)$-spheres $\mathcal{S}_{ij}$ pass through a point of this line if any of these spheres does.  
If such a point exists, we call it a \textit{generalized isodynamic point} of $\Sigma$ with respect to $P$. If there is exactly one generalized isodynamic point, then it is a point on the circumcircle. If there are two such points, they are inversives with respect to the circumcircle.\vspace*{1mm}

\begin{figure}[!t]
\includegraphics[height=8cm]{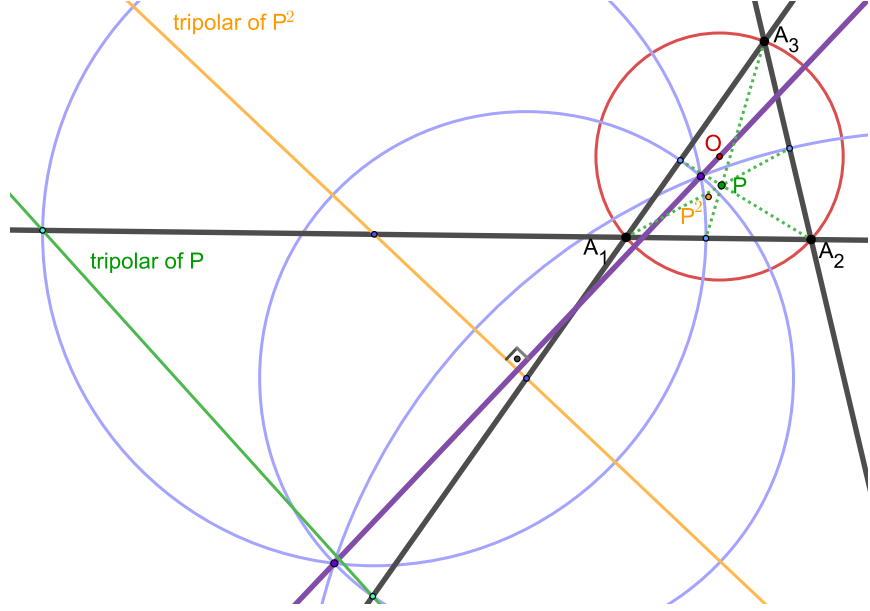}
\caption{Three generalized Appollonian circles. All figures were created with the software program GeoGebra [25].}
\end{figure}

The "classical" case: $P = \textrm{incenter}\, I$. In this case, the spheres $\mathcal{S}_{ij}$ are the (proper) \textit{Apollonian spheres}. $\check{\textrm{S}}\textrm{ruba}\check{\textrm{r}}$'s paper [21] suggests that the line through $O$ perpendicular to the $\Sigma$-polar plane of $K = I^2$ always meets an Apollonian sphere at two points, which are then called \textit{isodynamic points}. However, this is not the case, not even for $n=3$. We present a counter example. There exists a tetrahedron with sidelengths $(d_{12},d_{13},d_{14},d_{23},d_{24},d_{34}) = (13,11,9,12,5,11)$. We will show that the Apollonian spheres $\mathcal{S}_{12}, \mathcal{S}_{13}, \mathcal{S}_{23}$ of this tetrahedron do not meet. Let $a_i$ denote the area of the facet opposite vertex $A_i$; then  $a_1=6\sqrt{21}, a_2=\frac{9}{4}\sqrt{403}, a_3=\frac{9}{4}\sqrt{51}, a_4 = 6\sqrt{105}$. The line $A_4 \sqcup K$ meets the sideplane $A_1 \sqcup A_2 \sqcup A_3$ at the point $R = [a_1^{\;2}{:}a_2^{\;2}{:}a_3^{\;2}{:}0]$. The  intersection of the spheres $\mathcal{S}_{12}, \mathcal{S}_{13}, \mathcal{S}_{23}$ with the plane $A_1 \sqcup A_2 \sqcup A_3$ are the generalized Apollonian circles of the point $R$ and triangle $\Delta = A_1A_2A_3$. P. Yiu [23] gives a criterion to decide whether or not these circles have points in common: They have no common points precisely when the point \vspace*{1mm}\\
\centerline{$\displaystyle Q = \big{[}\, d_{23}^{\;\,2}\big{(}\frac{d_{23}^{\;\,2}}{a_1^{\;2}}-\frac{d_{13}^{\;\,2}}{a_2^{\;2}}- \frac{d_{12}^{\;\,2}}{a_3^{\;2}} \big{)}:d_{13}^{\;\,2}\big{(}\frac{d_{13}^{\;\,2}}{a_2^{\;2}}-\frac{d_{12}^{\;\,2}}{a_3^{\;2}}- \frac{d_{23}^{\;\,2}}{a_1^{\;2}} \big{)} : d_{12}^{\;\,2}\big{(}\frac{d_{12}^{\;\,2}}{a_3^{\;2}}-\frac{d_{23}^{\;\,2}}{a_1^{\;2}}- \frac{d_{13}^{\;\,2}}{a_2^{\;2}} \big{)} : 0\,\big{]}$}\vspace*{1mm}\\
lies outside the circumcircle of $\Delta$.\\
Since the distance between $Q = [\frac{3326952}{4504043},\frac{25180529}{27024258},-\frac{18117983}{27024258},0]$ and the circumcenter of $\Delta$, $O_\Delta = [\frac{73}{210},\frac{121}{315},\frac{169}{630},0]\,$, is bigger than the circumradius of $\Delta$, the Apollonian spheres $\mathcal{S}_{12}^{\,}, \mathcal{S}_{13}, \mathcal{S}_{23}$ of$^{\,}$ this$^{\,}$ tetrahedron do not meet, and isodynamic points, when defined as common points of the Apollonian spheres, do not exist for this tetrahedron.
\vspace*{2mm}\\
\textit{Remark}. C. Pohoata and V. Zajic [20] presented a generalization of the Apollonian circles different from that introduced by P. Yiu.  \\

\noindent\textbf{4. Isogonic points of a simplex.}\hspace*{\fill}\vspace*{1 mm} \\
From [4] we adopt the following terminology: A simplex is \textit{regular} if all its edges are the same length, a simplex is \textit{equiareal} if all its facets have the same $(n{-}1)$-volume, and it is \textit{equifacetal} if all its facets are congruent.\\
In dimension $3$, every equiareal simplex is equifacetal; but in higher dimensions this is not the case.\\
Let us define:\vspace*{0.5mm}\\
\centerline{A point $P$ is an \textit{isogonic point} of an $n$-simplex $\Sigma$ if its antipedal simplex is  equiareal.}\\
It follows immediately, that a point $P$ is an isogonic point of an $n$-simplex $\Sigma$ if and only if its mirror image $\Sigma^\star$ in an $(n{-}1)$-sphere with center $P$ is equiareal. The simplex $\Sigma^\star$ and the antipedal simplex are similar simplices. \vspace*{1.5mm}\\
Before moving to higher dimensions $n$, we shall see that our definition is adequate for dimension $n=3$: If the antipedal simplex is equarial, it is also equifacetal, and the four triads of lines $(P\sqcup A_i) \cup (P\sqcup A_j) \cup (P\sqcup A_k)$, $1\le i<j<k\le 4$, are congruent in pairs. On the other hand, this congruence also implies that the antipedal simplex is equifacetal. But it is not the congruence, it is the equiangularity that justifies the name \textit{isogonic point}. Moreover, if $P$ is an isogonic point inside the tetrahedron $\Sigma$, then this point is the Fermat-Torricelli point of $\Sigma$, see [1]. Let us call the Fermat-Torricelli point $F$. \vspace*{-2mm}\\

In 1937 A. Weiszfeld published a method to calculate the barycentric coordinates of $F$. Starting from a point $F_0 = [f_{0,1}:\dots:f_{0,n{+}1}]$ with $f_{0,1}\cdot {\dots} \cdot f_{0,n{+}1}\ne 0$, $F$ is the limit of an iteration process with a step: $F_{i{+}1} = [f_{i{+}1,1}:\dots:f_{i{+}1,n{+}1}], f_{i{+}1,k} = f_{i,k}/d(F_{i},A_i)$, see [22,{\,}4,{\,}14,{\,}15].  A geometric description of this process is given in the next section.\\

\noindent\textbf{5. $Z^\star$-correspondence.}\hspace*{\fill}\vspace*{1 mm} \\
Let $P = [p_1{:}\cdots {:}p_{n+1}]$ be a point not on any $(n{-}1)$-sideplane of $\Sigma$, $\mathcal{S}$ an $(n{-}1)$-sphere with center $P$, and let $\Sigma^\star = A_1^{\,\star}A_2^{\,\star}\cdots A_{n{+1}}^{\;\,\star}$ be the polar simplex of $\Sigma$ with respect to $\mathcal{S}$. The simplices $\Sigma$ and $\Sigma^\star$ are orthologic,  and both orthologic centers coincide at $P$.\\ Therefore the barycentric coordinates of $P$ with respect to $\Sigma^\star$ agree with the barycentric coordinates of $P$ with respect to $\Sigma$.\\
\textit{Proof} of the last statement. We can assume that $p_1{+}\cdots {+}p_{n+1}=1$ and that $\mathcal{S}$ has radius 1. It suffices to show that $\sum_{i} p_i A_{\,i}^\star = P$. First observe that the vector $\displaystyle\mathbf{n}_i := \frac{\textrm{sgn}(p_i)}{d(P,A_{\,i}^\star)}(A_{\,i}^\star-P)$ is an outward unit normal vector of the $(n{-}1)$-dimensional surface of $\Sigma$. When we denote the $(n{-}1)$-volume of the facet opposite vertex $A_i$ by $a_i$, then $\sum_{\,i} a_i \mathbf{n}_i$ is the zero vector. It follows:\vspace*{1mm}
\[
\begin{split} \sum_{\,i} p_i A_{\,i}^\star &= \sum_{\,i}\, p_i P\;  +\;  \sum_{\,i}\,  p_i\,(A_{\,i}^\star-P)\\
&=\hspace*{4.5mm} P\hspace*{5mm} +\; \sum_{\,i}\,\Big{(}\,\frac{a_i}{n\, d(P,A_{\,i}^\star)\sum_{\,j}a_j}\Big{)}\Big{(}d(P,A_{\,i}^\star)\,\mathbf{n}_i)\Big{)}\\
&=\hspace*{4.5mm}P\hspace*{5mm} +\; \frac{1}{n\,\sum_{\,i}a_i}\,\sum_{\,i} a_i \mathbf{n}_i\;\; =\;\; P.
\end{split}
\]

Now we take a point $Z^\star$ having homogeneous coordinates $z_{\;1}^{\star}{:}\cdots{:}z_{\;n{+}1}^{\star}$ with respect to $\Sigma^\star\!$. If $Z^\star$ is different from $P$, 
its polar $(n{-}1)$-plane with respect to $\mathcal{S}$ is a hyperplane in $\mathcal{A}$ and is called the {$Z^\star\textit{-transversal}$} of $P$. The $\Sigma$-pole of this hyperplane is called the $Z^\star\!$-\textit{correspondent} of $P$, cf. [6]. It makes sense to choose the centroid $G$ of $\Sigma$ as the $P^\star\!$-{correspondent} of $P$. \\
Let us denote the $Z^\star\!$-correspondent of $P$ by $P\# Z^\star$, then \vspace*{1mm}\\
\centerline{$\displaystyle P\# Z^\star = \big{[}\frac{p_1}{z_{\;1}^{\star}}:\cdots:\frac{p_{n+1}}{z_{\;n{+}1}^{\star}}\big{]}\hspace*{15mm}(\star\!\star\!\star)\;\;$.}\vspace*{0mm}\\

\begin{figure}[!t]
\includegraphics[height=9cm]{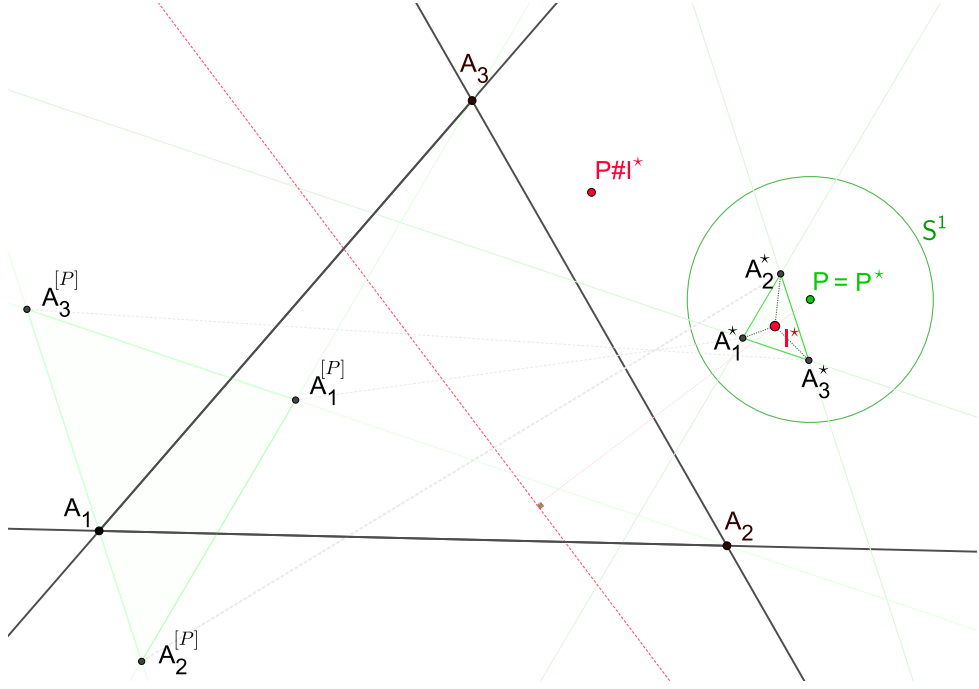}
\caption{The $I^\star$-transversal and the $I^\star$-correspondent of a point $P$. \hspace*{16mm} Triangle $A_1^{[P]}A_2^{[P]}A_3^{[P]}$ is the antipedal triangle of $P$.}
\end{figure}

\noindent \textit{Proof}. Here only an outline of a proof of this equation is presented; a more detailed proof is given in [6].\\
We assign to each point $R\in \mathcal{A}$ and each hyperplane $h$ not passing through $P$ a real number $\lambda_P =\lambda_P (R,h)$ by\vspace*{1mm}\\
$\hspace*{15mm}\lambda_P (R,h) = 0\,$, if $R = P$ or if the line $P\sqcup R$ does not meet $h$\vspace*{-1mm}\\
and \vspace*{-1mm}\\
$\hspace*{15mm}\lambda_P (R,h) = \frac{1}{t}\,$, if the point  $P+t\,(R-P)$ is the intersection of $P\sqcup R$ and $h$.\vspace*{1mm}\\
Some properties of $\lambda_P$ are listed below:\\
$\bullet$ For any points $R_1, R_2$, for any real numbers $t_1, t_2$ and for an hyperplane $h$ of $\mathcal(A)$,\\
\hspace*{25mm}  $\lambda_P(t_1 R_1+t_2 R_2, h) = t_1 \lambda_P(R_1,h) +t_2 \lambda_P(R_1,h)$. \\
$\bullet$ A point $R$ lies on a hyperplane $h$ if and only if $\lambda_P (R,h) = 1$. \\
$\bullet$ Let us denote the $(n{-}1)$-sideplane of $\Sigma$ opposite of the vertex $A_i$ by $h_i, i=1,\dots,n{+}1$.\\
\hspace*{2.5mm}If $P=[p_1,\dots,p_{n+1}]$, then  $\lambda_P (A_i,h_i) = 1-1/p_i$. It can be easily checked that the $i$-th\\\hspace*{2.5mm}barycentric coordinate of the point $P+\frac{p_i}{p_i-1}(A_i-P)$ is zero.\vspace*{-2.5mm}\\

If $\psi$ is the mapping that assigns to each point $R$ its $\mathcal{S}\textrm{-}$polar hyperplane, then \\
\hspace*{17mm}{$\psi(A_{\,i}^\star) = h_i = \{R\,|\,\lambda(R,h_i)=1\}$,}\\                                                             
\hspace*{17mm}$\psi(z_{\,i}^\star A_{\,i}^\star) = \{R\,|\,z_{\,i}^\star \lambda_P(R,h_i)=1\}$\\
\hspace*{8mm} and \hspace*{2mm}$\psi(\sum_{\,i}z_{\,i}^\star A_{\,i}^\star) = \{R\,|\,\sum_{\,i}z_{\,i}^\star \lambda_P(R,h_i)=1\}$.\vspace*{1mm}\\
\noindent\hspace*{5mm} We calculate the intersection of a sideline $A_i\sqcup A_j$ of $\Sigma$ with the hyperplane $\psi(Z^\star)$. For simplicity, we take $(i,j)=(1,2)$ and determine the real number $x$ such that\\ $\lambda_P(x A_1+ (1{-}x) A_2, \psi(Z^\star))=1$. \vspace*{0mm}
\[
\begin{split}
1 &= \lambda_P(x A_1+ (1{-}x) A_2, \psi(Z^\star))\\
&= x (\sum_k z_{\;k}^\star \lambda_P(A_1,h_k)) +(1{-}x)(\sum_k z_{\;k}^\star \lambda_P(A_2,h_k))\\
&= x (1-z_{\;1}^\star + z_{\;1}^\star(1-\frac{1}{p_1})) + (1{-}x)(1-z_{\;2}^\star + z_{\;2}^\star(1-\frac{1}{p_2}))\\
&= x(1-\frac{z_{\;1}^\star}{p_1})+(1-x)(1-\frac{z_{\;2}^\star}{p_2})\;.
\end{split}
\]
As a result we get $\displaystyle x= \frac{p_1/z_{\;1}^\star}{p_1/z_{\;1}^\star-p_2/z_{\;2}^\star}$ and $\displaystyle 1{-}x = -\frac{p_2/z_{\;2}^\star}{p_1/z_{\;1}^\star-p_2/z_{\;2}^\star}.\hspace*{20mm}\Box$\\

\begin{figure}[!t]
\includegraphics[height=9cm]{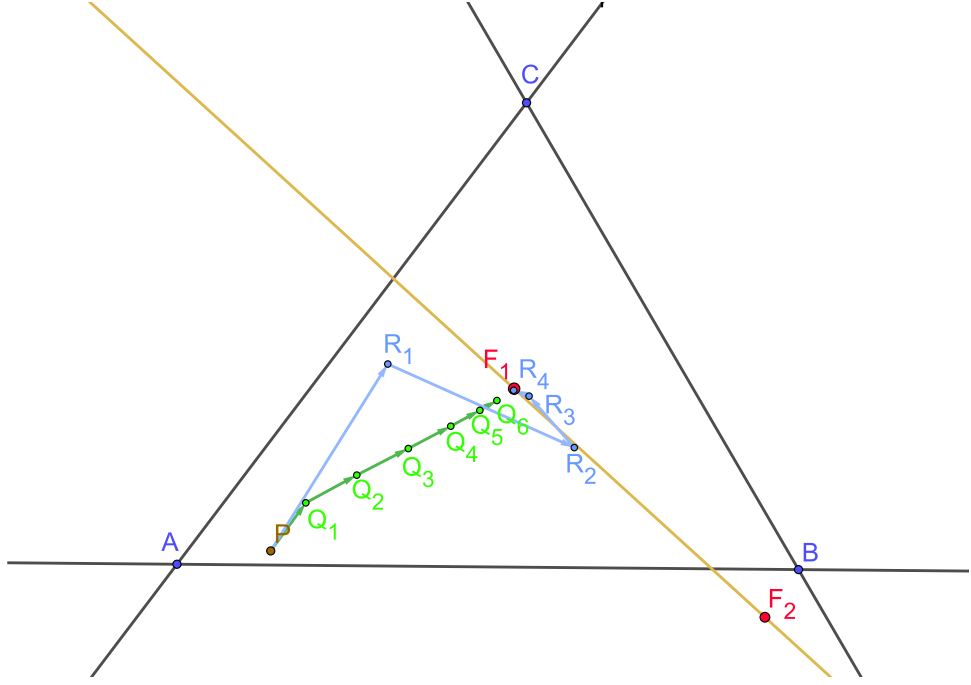}
\caption{The first points of the sequences $(Q_i)_{i\in\mathbb{N}}$ and $(R_i)_{i\in\mathbb{N}}$ with $Q_0 = R_0 = P$ and $Q_{i{+}1} = Q_i\#I^\star,  R_{i{+}1} = R_i\#K^\star$ are shown. $F_1$ and $F_2$ are the two isogonic points.}
\end{figure}
\vspace*{-1 mm}

In the following we are particularly interested in the $Z^\star$-correspondents for $Z^\star = G^\star$ (centroid of $\Sigma^\star$), $Z^\star = I^\star$ (incenter of $\Sigma^\star$) and $Q^\star = K^\star$ (symmedian of $\Sigma^\star$):\vspace*{1mm}\\
$\bullet\; P\# G^\star = P$\\
$\displaystyle \bullet\; P\# I^\star \,= [\frac{p_1}{a_{\,1}^\star}:\cdots:\frac{p_{n+1}}{a_{\,n+1}^\star}] = [\frac{\textrm{sgn}(p_1)}{d(P,A_1)}:\cdots:\frac{\textrm{sgn}(p_{n{+}1)}}{d(P,A_{n+1})}]$, $a_{\,i}^\star$ = $(n-1)$-volume of \hspace*{3mm}the facet of $\Sigma^\star$ opposite $A_{\,i}^\star$.\\
\hspace*{3mm}Explanation: $\displaystyle p_i:p_j = p_{\,i}^\star : p_{\,j}^\star = \textrm{sgn}(p_i)\frac{a_{\,i}^\star}{d(P,A_i)}: \textrm{sgn}(p_j)\frac{a_{\,j}^\star}{d(P,A_j)}.$\\
$\displaystyle \bullet\, P\# K^\star = [\frac{p_1}{(a_{\,1}^\star)^2}:\cdots:\frac{p_{n+1}}{(a_{\,n+1}^\star)^2}]=[\frac{1}{p_1(d(P;A_1))^2}:\cdots:\frac{1}{p_{n+1}(d(P;A_{n{+}1}))^2}]\,.$\\

For a point $P = [p_1,\dots,p_{n+1}]$ not on any sideplane of $\Sigma$ consider the sequences $(Q_i)_{i\in\mathbb{N}}$ and $(R_i)_{i\in\mathbb{N}}$ given by $Q_0 = R_0 = P$ and $\displaystyle Q_{i{+}1} = \big{[}\frac{1}{d(Q_i;A_1)}{:\dots:}\frac{1}{d(Q_i;A_{n+1})}\big{]},\\
\hspace*{48,5 mm}R_{i{+}1} = \big{[}\frac{1}{|p_1|(d(R_i;A_1))^2}:\cdots:\frac{1}{|p_{n+1}|(d(R_i;A_{n{+}1}))^2}\big{]}$.\vspace*{1 mm}\\
All points $Q_i, R_i, i>0,$ are points inside $\Sigma$, and $Q_{i+1} = Q_i\# I^\star, R_{i+1} = R_i\# K^\star$. For $Q_1$ and $R_1$ we have  $Q_1=P\# Z^\star$ with  $Z^\star = [\textrm{sgn}(p_1)\,a_{\,1}^\star:\cdots:\textrm{sgn}(p_{n+1})\,a_{\,{n+1}}^\star]$  and $R_1=P\# Z^\star$ with $Z^\star = [\textrm{sgn}(p_1)(a_{\,1}^\star)^2:\cdots:\textrm{sgn}(p_{n+1})(a_{\,{n+1}}^\star)^2]$. 
	Both sequences converge to the Fermat-Torricelli point $F$, even if this is a vertex of $\Sigma$. It seems that the second sequence in comparison to the first converges twice as fast (see Figure 3). But more importantly, the second sequence avoids calculating square roots. If the Fermat-Torricelli point of a simplex lies inside the simplex, then it is an isogonic point.\\

\noindent\textbf{6. A $3$-simplex with five isogonic points.}\hspace*{\fill} \vspace*{1 mm}\\  
In dimension $n>2$, simplices, in general, have more than two isogonic points. We present a $3$-simplex with five isogonic points. Consider in the Euclidean space $(\mathbb{R}^4, \cdot)$ the simplex $\Sigma$ with the four vertices\vspace*{1mm}\\
 \centerline{$A_1 = \left(  
\begin{array}{c}
0\\
0\\
0\\
\end{array}
\right) ,
 \,A_2 =\left(
\begin{array}{c}
6\\
0\\
0\\
\end{array}
\right), 
\,A_3=\left( \begin{array}{c}
0\\
8\\
0\\
\end{array}
\right), 
\, A_4 = \left(
\begin{array}{c}
2\\
2\\
6\\
\end{array}
\right)\,$.}\\ 

In order to calculate the isogonic points, we first calculate their isogonally conjugated points. These are the points with an equiareal pedal simplex. For a point $P$ which is not a vertex of $\Sigma\,$ let $\Sigma_P$ denote the pedal simplex of $P$ and let $G_P$ and $I_P$ denote the centroid and the incenter of $\Sigma_P$, respectively. We define a sequence of points $(P_i\,)_{\,i\in \mathbb{N}}$ by $P_0 := P$ and $P_{n+1}:=P_{n}+G_{\!P_n}-I_{\!P_n}\,$. Obviously, this sequence is a constant sequence if and only if $P$ is the isogonal conjugate of an isogonic point. Experiments with the interactive geometry tool GeoGebra [25] suggest that each of these sequences converges to one of these limit points 
(rounded to $12$ decimal places):\vspace*{-2mm}\\

\centerline{$L_0 =[ 0.266996565955 , 0.275481800939 , 0.217355830792 , 0.240165802314]$,}
\centerline{$\;\;\,L_1=[-4.180629474014 , 2.569387212447 ,1.602113038329 , 1.009129223238]$,}
\centerline{$\;\;\,L_2=[1.193250865914, -1.252645952150, 0.354761022780, 0.704634063455]$,}
\centerline{$\;\;\,L_3=[0.713260932730, 0.358215195120, -0.616627271982, 0.545151144132]$,}
\centerline{$\;\;\,L_4=[0.657546390333, 0.802131717931, 0.639088262811, -1.098766371077]$.}\vspace*{2mm}

\noindent The pedal triangles of these points have facets with an area\vspace*{1mm}\\
\centerline{$a_0= 2.404772767371,\, a_1=122.125536031480,\, a_2=19.392997370805,$}
\hspace*{11.2mm}$a_3= 9.848601171111,\, a_4=18.965046082427\,$.\vspace*{-1mm}\\

\noindent The isogonic points of $\Sigma$ are\vspace*{-2mm}\\

\centerline{$F_0 =[0.369979160947, 0.229493293826, 0.163611619856, 0.236915925371]$,}
\centerline{\;\;\,$F_1 =[-0.297000489955, 0.309278164652, 0.279002561033, 0.708719764270]$,}
\centerline{\;\;\,$F_2 =[0.388102931405, -0.236608485604, 0.469943106828, 0.378562447371]$,}
\centerline{\;\;\,$F_3 =[0.382915343108, 0.487963317698, -0.159452369671, 0.288573708865]$,}
\centerline{\;\;\,$F_4 =[0.645021938255, 0.338403751068, 0.238914519123, -0.222340208446]$.}\vspace*{2mm}

\noindent The antipedal triangles of these points have facets with an area\vspace*{1mm}\\
\centerline{$\tilde{a}_0= 241.637142362610,\, \tilde{a}_1=60.087819904352,\, \tilde{a}_2=31.387257487815,$}
\hspace*{10.4mm}$\tilde{a}_3= 5.647726265255,\hspace*{4mm} \tilde{a}_4=31.003305976553\,$.\vspace*{-1mm}\\

\noindent Finally, we present the two isodynamic points of $\Sigma\,$:\vspace*{-2mm}\\

\centerline{$J_1 = [0.206439675828, 0.327649375007, 0.263085414624, 0.20282553454],\,\;\,$}
\centerline{$\,\;\,\,\;\,J_2 = [2.954833710960, -0.575606610593, -1.403778427224, 0.024551326857]$.}

\end{document}